\documentclass[]{amsart}

\theoremstyle{plain}
\newtheorem*{main theorem}{Main Theorem}
\newtheorem*{theorem*}{Theorem}
\newtheorem*{emptytheorem*}{}
\newtheorem{theorem}{Theorem}

\newtheorem*{proposition*}{Proposition}
\newtheorem{lemma}{Lemma}

\newtheorem{corollary}{Corollary}

\theoremstyle{definition}
\newtheorem{definition}{Definition}
\newtheorem*{definition*}{Definition}

\theoremstyle{remark}

\newtheorem*{remark*}{Remark}

\begin{document}
\title[Hyperbolicity estimates for random maps]
{Uniform hyperbolicity for random maps \\
with positive Lyapunov
exponents}

\author{Yongluo Cao}
\address{Dept. of Mathematics, Suzhou University,
Suzhou 215006, Jiangsu, P.R. China}
 \email{ylcao@suda.edu.cn}

\author{Stefano Luzzatto}
\address{Dept.  of Mathematics, Imperial College, 180 Queen's
Gate, London SW7 2AZ, UK}
\email{Stefano.Luzzatto@imperial.ac.uk}
 \urladdr{http://www.ma.ic.ac.uk/\~{}luzzatto }

\author{Isabel Rios}
\address{Universidade Federal Fluminense, Niteroi, RJ, Brazil.}
\email{rios@mat.uff.br}

\thanks{IR was partially supported by CAPES and FAPERJ (Brazil).
YC was partially supported by  NSFC(10571130), NCET, 
and SRFDP of China
and the Royal Society. SL was partially supported by EPSRC grant 
GRT0969901.
The authors would like thank M. Benedicks and 
M. Viana for their suggestions and encouragement.}

\date{31 March 2007. Minor revisions 2 August 2007.}


\begin{abstract}
We consider some general classes of random dynamical systems and show 
that a priori very weak nonuniform hyperbolicity conditions
actually imply uniform hyperbolicity.
\end{abstract}    

\maketitle

\section{Introduction and statement of results}

In this paper we consider smooth random
dynamical systems $F$ over an abstract dynamical systems
$(\Omega, \mathcal{F},\mathbb{P}, \theta)$, where $(\Omega,
\mathcal{F},\mathbb{P})$ is a complete probability
space and $\theta: \Omega \to \Omega $ 
is a $\mathbb{P}$ preserving ergodic invertible
transformation. More specifically, we have a 
skew-product
$$
F : \Omega \times M \to \Omega \times M,
$$
given by $$F(\omega,x) = (\theta(\omega), \phi_{\omega}(x))$$
where   $M$ is a compact manifold endowed with a Riemannian metric
 which induces a norm $|\cdot |$ on the tangent space and  a volume
 form that we call
 Lebesgue measure. Throughout the paper we
 suppose that for $\mathbb{P}-a.e $
 \( \omega \),
 \[
 \phi_{\omega}: M \to M 
 \]
 is a $C^1$  local diffeomorphism. We let 
\[
|D\phi_{\omega}|=\sup_{x \in M} |D_x \phi_{\omega}|\quad
\text{and}
\quad 
|D\phi^{-1}_{\omega}|=\sup_{x \in M} |D_x \phi^{-1}_{\omega}|
\]
and assume standard integrability conditions
  \begin{equation}\label{intass}
  \int\limits_{\Omega} |D\phi_{\omega}|d\mathbb{P} <
  \infty\quad\text{and}\quad 
  \int\limits_{\Omega}
  |D\phi_{\omega}^{-1}|d\mathbb{P} < \infty
  \end{equation}
Notice that these conditions are not automatic since we do not assume 
that \( \phi_{\omega} \) depends continuously on \( \omega \) in any
way.

 Random maps of this kind have been extensively studied from various
 points of view, such as the existence and
 properties of invariant measures and equilibrium states 
 \cite{KhaKif96, ArbMatOli03} and the
continuity properties of the entropy \cite{LiuZha03},
 see \cite{Kif98, Liu01} for an extensive survey and
 references. 
Many results depend 
on some \emph{hyperbolicity} of the random maps
or, in the language of skew-products, on some \emph{hyperbolicity
in the fibres}. Our main goal
in this paper is to show that in several cases, 
\emph{uniform} hyperbolicity estimates
can be obtained from some a priori strictly weaker
\emph{nonuniform hyperbolicity} assumptions.

\subsection{Basic definitions}

\subsubsection{Random continuous functions}
A function $ f : \Omega \times M  \to R $ is a random continuous
function if 
\begin{enumerate}
    \item 
    $ x \mapsto f(\omega , x) $ is a  continuous
function for a.e. $\omega \in \Omega$;
\item $\omega
\mapsto f(\omega, x ) $ is  measurable for all $  x \in M $; 
\item
$\omega \mapsto \sup_{x \in M }|f(\omega,x)| $ is integrable with
respect to $\mathbb{P}$. 
\end{enumerate}

\subsubsection{Topology on the space of measures}
We let $ \mathcal M_{\mathbb{P}} (F)$ denote   all $ F $-invariant
probability measures
 on $\Omega\times M $ whose marginal on $ \Omega $ coincide with
  $ \mathbb{P}$ (such measures can be characterized in term of
  their disintegrations $\mu_{\omega}$ by
  $\phi_{\omega} (\mu_{\omega}) = \mu_{\theta \omega} a.s.$).
We equip $ \mathcal M_{\mathbb{P}} (F)$ with the
smallest topology such that
$$ \nu \to \int_{\Omega} \int_M f(\omega, x) 
d\mu_{\omega}(x)d\mathbb{P}(\omega)
=\int_{\Omega \times M }f(\omega, x )d\mu(\omega, x) $$ 
is continuous for every random continuous function $f$.
We let
 \[
 \mathcal E_{\mathbb{P}}(F) \subset \mathcal M_{\mathbb{P}}(F)
 \]
 denote
the subset of ergodic measures.

\subsubsection{Fibrewise Lyapunov exponents}

For $\omega \in \Omega$, let $\phi^{(0)}_{\omega}$ be
the identity map on $M$ and, for $k \in \mathbb N $, define
$\phi^{(k)}_{\omega}$ by
$$
\phi_{\omega}^{(k+1)} = \phi_{\theta^k(\omega)} \circ
\phi_{\omega}^{(k)}.
$$
Then we can define a family of iterates of \( F \) by
\[
F^{n}(\omega,x) = (\theta^{n}(\omega), \phi_{\omega}^{(n)}(x)).
\]
The derivative map of $\phi$ along the $M$ direction gives a
cocycle 
\[
(\omega,x,k) \to D_x \phi_{\omega}^{(k)}
\]
from 
$\Omega \times M \times \mathbb N $ to $GL(m,R)$ where $m
= \dim M$.   
\begin{definition}
For each $\omega \in \Omega, x \in M $ and $v \in T_xM $, we say
that
 $$
 \lambda( \omega ,x,  v ) = \lim_{k \to \infty } \frac1k \log ( |
 D_x\phi_{\omega}^{(k)} (v) | ), 
 $$
if the limit exists,  is the \emph{fibrewise Lyapunov exponent}
 associated to the point \( (\omega,x)  \) and the vector \( v \).
 \end{definition}
By Oseledec's theorem (see \cite{Ose68, Rue82}) the limit exists
for $\nu$-almost all
 $(\omega,x)$ for any   F-invariant
probability measure \( \nu \) and therefore 
for a.e. \( (\omega, x) \) there are 
 real numbers  $\lambda_1(\omega,x) \le \lambda_2(\omega,x) \le
 \cdots \le \lambda_m(\omega,x)$ which 
 are the fibrewise Lyapunov exponents corresponding to different
directions in \( T_{x}M \). If \( \nu \) is ergodic, these numbers
 are constant $\nu$-almost everywhere
 and we denote them as $\lambda_1(\nu) \le
 \cdots \le \lambda_m(\nu)$.
 
 \subsection{Random expanding maps}
 
 \subsubsection{Positive Lyapunov exponents}
 In this paper we shall be particularly interested in the case in which
 the fibrewise Lyapunov exponents are positive. 
  Notice that in this case the definition implies that for all \(
  \varepsilon>0 \) sufficiently small, there exists a
  constant \( C(\varepsilon, \omega, x)>0 \) such that 
  \[ 
  |D_x\phi_{\omega}^{(k)} (v) | \geq C(\varepsilon, \omega, x)
  e^{(\lambda(\omega,
  x, v)-\varepsilon)  n} |v|
  \]
  for all \( n\geq 1 \). 
In particular, for an ergodic \( F \)-invariant measure \( \nu \) with all
  fibrewise Lyapunov exponents positive: 
  \(
  \lambda_m(\nu) \ge
 \cdots \ge \lambda_1(\nu)=:\lambda (\nu)> 0.
 \)
 This implies that, for \( \nu \) 
 almost all \( (\omega, x) \) and for all sufficiently small \(
 \varepsilon>0 \), there exists a constant \( C(\varepsilon, \omega,
 x)>0 \) such that 
 \[ 
 |D_x\phi_{\omega}^{(k)} (v) | \geq C(\varepsilon, \omega, x)
 e^{(\lambda-\varepsilon) n}|v|
 \]
 
 \subsubsection{Random uniform expansion}
In certain cases such an expansion estimate 
actually extends to all of \( M \) with a constant \( C \) independent
of the point \( x \).

\begin{definition}\label{expanding} A random map $F$ is called 
    \emph{random uniformly expanding} if
there exists a constant \( \lambda>0 \) and a 
tempered random variable $
C(\omega) > 0$ such that  for 
$\mathbb{P}$-\emph{almost all} $\omega \in \Omega $ and 
\emph{every} \( x\in
M\) we have
$$\|D_x\phi_{\omega}^n( v )\| \ge C(\omega) e^{\lambda n}\|v\|.$$
\end{definition}

Notice that this extends the usual definition of a \emph{uniformly
expanding} map to the random setting by requiring the expansion rate \( 
\lambda \) to 
be uniform in both \( \omega \) and \( x \) though still allowing the
constant \( C \) to depend (in a controlled way) 
on \( \omega \).
We recall that a
random variable $g : \Omega \to \mathbb R^{+}$  is \emph{tempered}
if 
\[
\lim\limits_{n \to \infty} \frac1n \log
g(\theta^n(\omega)) = 0 , \mathbb{P}-a.s.
\]
Our first result says that such a uniform expansion property
actually follows from an a priori weaker assumption.

 \begin{theorem}\label{th1}
Let \( F \) be a random map and suppose that all fibrewise Lyapunov
exponents are positive for all measures
\( \nu \in \mathcal E_{\mathbb{P}}(F) \)  . 
Then \( F \) is random uniformly expanding.
\end{theorem}

We emphasize that in our case the Lyapunov exponents are not
assumed to be uniformly bounded away from 0. Thus, a priori, we only
have that for every \( (\omega, x) \) in 
a subset of \( \Omega \times M \) of \emph{full
probability}, i.e. in a set which has full measure for every invariant
probability measure, there are constants \( C(\omega, x)>0 \) and \(
\lambda (\omega, x)>0 \) such that 
\( |D_x\phi_{\omega}^{(k)} (v) | 
\geq C(\omega, x) e^{\lambda (\omega, x) n} |v| \). 
Theorem~\ref{th1} says that the  expansion estimates 
actually hold for
every \( x \) and for constants \( C, \lambda\) independent of \( x \),
and thus in particular that all fibrewise Lyapunov exponents are
 uniformly bounded away from 0.

 \subsubsection{Deterministic case}
We remark that the results are non-trivial even in the special
case in which the \( \theta \)-invariant measure \( \mathbb{P}\)
is a Dirac-\( \delta \) measure supported on a single fixed point
\( \{p\} \). The setting stated above then reduces to the case in
which \( F: M \to M \) is a standard deterministic dynamical
system and an analogous result has been proved in
\cite{AlvAraSau03, Cao03,
 CaoLuzRioHyp}. The theorem we prove here represents a
significant generalization of these results and is obtained by a
different argument. The general question of the
uniformity of in principle nonuniform functions has also been
addressed in various contexts in other papers such as
\cite{StaStu00, Slo97}.

\subsubsection{Uniform bounds for expansion rates}
As an immediate corollary of Theorem~\ref{th1} we get the following
statement. 

\begin{corollary} Let \( F \) be a random map and suppose that 
    there exists
    tempered random variables $
    C(\omega) > 0$ and \( \lambda(\omega) \)  with
    $\int \log \lambda d\mathbb{P} > 0 $, 
    such that  for 
    $\mathbb{P} $ almost all $\omega \in \Omega $ and every \( x\in M \)
    we have
    $$\|D_x\phi_{\omega}^n( v )\| \ge C(\omega) \lambda^{(n)} 
    (\omega) \|v\|$$
    where $\lambda^{(n)}(\omega) = \lambda(\omega) \cdots
    \lambda(\theta^{n-1} \omega)$. 
Then \( F \) is random expanding. 
In particular $\lambda(\omega)>1 $
can be chosen constant.
\end{corollary}
\begin{proof} The assumption that 
    \( \int \log \lambda d\mathbb{P} > 0 \)
    implies   then for all measures
\( \nu \in \mathcal E_{\mathbb{P}}(F) \) all fibrewise Lyapunov
exponents are positive. Then Theorem~\ref{th1} implies the result.
\end{proof} 

\subsection{Random hyperbolic maps}
We now state versions of these results for the cases in which 
$ \theta $ is an invertible
transformation and 
$\phi_{\omega}$ is a  $C^1$ diffeomorphism for a.e. \( \omega \). 

\subsubsection{Random compact sets}

\begin{definition}
$\Lambda=\{\Lambda(\omega): \omega \in \Omega \}$ is a
\emph{random compact set}
if 
\begin{enumerate}
 \item 
$\Lambda(\omega) \subset M$ is compact for a.e. \( \omega \); 
\item 
$(x, \omega ) \to d(x,\Lambda(\omega))$ is measurable.
\end{enumerate}
Here $d$ is the Hausdorff distance on $M$.
A random compact nonempty set
$\Lambda=\{\Lambda(\omega): \omega \in \Omega \}$ is
\emph{invariant} under $F$ if 
\[
\phi_{\omega} \Lambda(\omega) =
\Lambda(\theta \omega)
\]
for a.e. $\omega \in
\Omega$.
\end{definition}

\subsubsection{Random uniform hyperbolicity}
\begin{definition}
A random, compact, $F$-invariant, nonempty, set
$\Lambda$ has a 
\emph{uniform tangent bundle splitting}  
if there exist 
\textit{i)} an open set $V$ with
a compact closure $ \bar{V}$, \textit{ii)}
a tempered random variable $\alpha > 0$ 
with $\int \log \alpha d\mathbb P <\infty$, 
and  \emph{iii)} subbundles $\Gamma^{1}(\omega)$ and $
\Gamma^{2} (\omega)$ of the tangent bundle $T\Lambda(\omega)$,
depending measurably on $\omega$  and continuously  on $x$, 
such that
\begin{enumerate}
\item There exist a measurable family of open set $U(\omega)$ such that
\begin{enumerate}
    \item 
    \(
\{x:d(x,\Lambda(\omega)) < \alpha(\omega)\} \subset U(\omega)
\subset V\);
\item
\( \phi_{\omega} U(\omega) \subset V\);
\item 
$\phi_{\omega}$ restricted to $U(\omega)$ is a diffeomorphism;
\item 
$\int\log^+\sup_{x \in U (\omega)}|D_x\phi_{\omega}|d\mathbb P <\infty$ and
$\int\log^+\sup_{x \in U (\omega)}|D_x\phi^{-1}_{\omega}|d\mathbb
P<\infty$.
\end{enumerate}
\item
\begin{enumerate}
    \item 
\( T\Lambda(\omega)=\Gamma^1(\omega)\bigoplus\Gamma^2(\omega)\);
\item 
\(D\phi_{\omega}\Gamma^1(\omega)=\Gamma^1(\theta
\omega)\) and \( D\phi_{\omega}\Gamma^2(\omega)=\Gamma^2(\theta
\omega)\);
\item 
\(\angle(\Gamma^1(\omega),\Gamma^2(\omega)) \ge
\alpha(\omega), \) for a.e \(\omega \).
\end{enumerate}
\end{enumerate}
\end{definition}

\begin{definition}[\cite{GKifer}] A random,  compact,  $F$-invariant, nonempty set
$\Lambda$  is a
\emph{random uniformly hyperbolic set} if it has 
a uniform tangent bundle splitting and 
there exists a constant \( \lambda > 0 \) and a
 tempered random 
variable $C>0$ such that 
for a.e $\omega$ and every $ n \in \mathbb{N}$ we have
$$
|D\phi_{\omega}^{(n)} \xi | \le C(\omega) e^{\lambda n}
|\xi| \quad \text{for} \quad \xi \in \Gamma^1(\omega)
\text{ and }
|D\phi_{\omega}^{(-n)} \eta | \le C(\omega)
e^{\lambda n}|\eta| \quad \text{for} \quad
\eta \in \Gamma^2(\omega).
$$
\end{definition}

\begin{theorem}\label{th2}
Let  $\Lambda$ be a
random,  compact,  $F$-invariant, nonempty, set with a uniform tangent
bundle splitting, and 
suppose  that for all measures \(
\nu \in \mathcal E_{\mathbb{P}}(F) \), all fibrewise Lyapunov
exponents restricted to $\Gamma^{1}$ are negative and all
fibrewise Lyapunov exponents restricted to $\Gamma^{2}$  are
positive. Then $\Lambda$ is a
random uniformly hyperbolic set for $F$.
\end{theorem}

Once again, we emphasize that this result is about showing that 
non-zero Lypaunov exponents on a full probability set 
actually imply uniform hyperbolicity
and thus, in particular, that all Lyapunov exponents are actually
uniformly bounded away from zero. 

Theorem~\ref{th2} follows immediately from Theorem~\ref{th1}
applied to each of the subbundles independently. In what follows we
therefore assume the setup and assumptions of Theorem~\ref{th1}.

 \section{Invariant measures on the unit tangent bundle}

Let $SM = \{(x, v ) \in TM : |v| = 1 \} $ denote the unit tangent
bundle over \( M \) and define the induced skew-product  tangent
map
\[
\widehat{TF} : \Omega \times  SM \to \Omega \times  SM
\]
by
$$
\widehat{TF} (\omega ,x, v ) = \left(\theta(\omega)
,\phi_{\omega}(x), \frac{D_x\phi_{\omega}(v)}
{|D_x\phi_{\omega}(v)|}\right).
$$
Since $\phi_{\omega}$ is  a $C^1$ local diffeomorphism, the
denominator in the definition  above never vanishes and hence this
map is well defined for all $(\omega,x, v ) \in \Omega \times SM
$. Extending the notation introduced above, we let $Pr(SM)$ denote
all probability measures  supported on $SM$ and  $ \mathcal
M_{\mathbb{P}} (\widehat{TF})$ denote all $ \widehat{TF}
$-invariant probability measures
 on $\Omega\times SM $ whose marginal on $ \Omega $ coincide with
  $ \mathbb{P}$ and  let
 $\mathcal E_{\mathbb{P}}(\widehat{TF}) \subset \mathcal M_{\mathbb{P}}(\widehat{TF})$  denote
the subsets of ergodic measures.
Since $ SM $ is compact, $ \mathcal M_{\mathbb{P}}(\widehat{TF})$ 
is compact in
the weak-star topology. Let
\[
\pi : \Omega \times  SM \to \Omega \times  M
\]
be the projection onto $ \Omega \times M $. 
We have $\pi
\circ \widehat{TF} = F \circ \pi $, and so if
 $m \in \mathcal M_{\mathbb{P}}(\widehat{TF})
 $, then $\pi^{*} m  = m \circ \pi^{-1} \in \mathcal M_{\mathbb{P}}(F).$
 Thus $\pi^{*} $
  defines a map
  \[
  \pi^{*}: \mathcal M_{\mathbb{P}}(\widehat{TF}) \to
  \mathcal M_{\mathbb{P}}(F).
  \]

\begin{lemma}\label{ergodic}
\(
\pi^{*}(\mathcal E_{\mathbb{P}}(\widehat{TF})) \subset  
\mathcal E_{\mathbb{P}}(F). \)
  \end{lemma}

 \begin{proof}
Let $A$ be  a measurable set  which is $F$ invariant.
 Then $\pi^{-1} A $ is a $\widehat{TF}$ invariant set. Since $m$ is
 ergodic,   $m(\pi^{-1} A ) = 0\)  or 1.
 Thus $ \nu (A) = \pi^{*} m (A) =m(\pi^{-1} A) = 0 \) or
 1. So  $\nu$ is ergodic.
\end{proof}

\section{Uniformly positive Lyapunov exponents}

We define the random  continuous function $\Phi
:\Omega \times SM \to R $ by
$$\Phi (\omega,x, v) = \log |D_x\phi_{\omega}(v)|.$$

\begin{lemma}\label{noname}
    There exists a measure 
    \( m^{*}\in  \mathcal M_{\mathbb{P}}(\widehat{TF})\)
    such that 
    $$
    \min_{m \in \mathcal M_{\mathbb{P}}(\widehat{TF})} \int_{\Omega
    \times SM} \!\!\!\!\!\!\!\!  \Phi   dm  = 
    \int_{\Omega \times SM} \!\!\!\!\!\!\!\! 
    \Phi  dm^* =: \Lambda > 0.
    $$
 In particular, all the fibrewise Lyapunov exponents of all invariant measures
 are uniformly bounded away from 0.    
\end{lemma}

\begin{proof}
The existence of a minimizing measure \( m^{*} \) follows immediately 
from the fact that   $\mathcal M_{\mathbb{P}}(\widehat{TF})$ is compact 
and by noticing that
$\Phi$ is  a random continuous function
    on $\Omega \times SM$ and therefore $\int \Phi dm  $ is continuous
    function on $\mathcal M_{\mathbb{P}}(\widehat{TF})$. 
    Therefore it 
only remains to show that \( \int \Phi dm^{*}>0 \) or,
equivalently, that \( \int\Phi dm > 0 \) for any 
 \( m\in \mathcal M_{\mathbb{P}}(\widehat{TF}) \). 
Moreover, 
by the Ergodic Decomposition Theorem we can assume without loss of
generality that \( m \) is ergodic. 

Thus let  $m \in \mathcal E_{\mathbb{P}}(\widehat{TF})$ and 
$\nu = \pi^{*} m \in \mathcal E_{\mathbb{P}}(F)$. Notice
that \( \pi \) maps full measure sets for \( m \) to full measure sets
for \( \nu \). 
By Birkhoff's Ergodic Theorem we have , for \( m \) almost every \(
(\omega, x, v) \), 
\[ 
\int_{\Omega \times SM} \!\!\!\!\!\! \!\!\!\!\!\! 
\Phi(\omega,x,v) dm = \lim_{n \to \infty}
 \frac1n \sum_{i=0}^{n-1} \Phi ((\widehat{TF})^i(\omega,x,v)) 
 \]
By the definition of \( \Phi \) we have 
\[
    \sum_{i=0}^{n-1} \Phi ((\widehat{TF})^i(\omega,x,v)) =
\log |D_x\phi_{\omega}^{(n)}(v)|.
\]
and therefore
\[ 
\lim_{n \to \infty}
 \frac1n \sum_{i=0}^{n-1} \Phi ((\widehat{TF})^i(\omega,x,v)) 
 = 
 \lim_{n \to \infty } \frac1n \log |D_x \phi_{\omega}^{(n)}v|.
\]
Applying Birkhoff's Theorem again, the limit on the right converges 
to \( \lambda (\omega,x,v)  \) which is \( >0 \) by our assumptions 
that all
fibrewise Lyapunov exponents are positive. 
\end{proof}

\section{Uniform hyperbolicity}

In the previous section we showed that all Lyapunov exponents are
uniformly bounded away from zero. We now need to extend the
corresponding expansion estimates to every point \( x \in M \). 

\begin{lemma}
For any \( \Lambda > \lambda > 0 \) we have that for a.e. \( 
\omega \) there exists a constant \( C(\omega)>0 \) such that for all \( 
x\in M \), \( v\in T_{x}M \), and \( n\geq 1 \)
\[ 
|D_{x}\Phi^{(n)}_{\omega}v|\geq C(\omega) e^{\lambda n}|v|.
\]
\end{lemma}
Notice that this is not quite the end result since we still need to
prove that \( C(\omega) \) is tempered. We shall do this in the next
section. 
\begin{proof}
The statement follows immediately from the fact that 
for a.e. \( \omega \) we have 
\[ 
\lim_{n\to\infty}\frac{1}{n}\min_{(x,v) \in SM} \{\log
|D_x\phi^{(n)}_{\omega}v|\} = \Lambda.
\]
To prove this, we show first of all that the limit exists and
is independent of \( \omega \), then we show that it is equal to \(
\Lambda \). 
\subsubsection*{Existence of the limit}
To get the existence of the limit, let 
 \[
 A_n(\omega) = \min\limits_{(x,v) \in SM} \log
|D_x\phi^{(n)}_{\omega}v|.
\]
Then $ A_{n+m}(\omega) \geq
A_n(\omega) + A_m(\theta^n \omega)$. Therefore the sequence \(
\{A_{n}\} \) is supadditive, the sequence \( \{-A_{n}\} \) is
subadditive and,  from the  subadditive ergodic
theorem \cite{King} and the ergodicity of $\mathbb{P}$ there exists a 
constant \( A \) such that 
\[
\lim_{n\to\infty}\frac{1}{n}\min_{(x,v) \in SM} \{\log
|D_x\phi^{(n)}_{\omega}v|\} = A 
\]
for a.e. \( \omega \). 

\subsubsection*{Upper bound}
From the previous section we know that 
\( n^{-1}\log |D_x\phi^{(n)}_{\omega}v| \) converges to \( \Lambda \) 
for some points (indeed, a set of points of full measure for the
minimizing measure \( m^{*} \)) and therefore we must have \( A\leq
\Lambda \). 

\subsubsection*{Lower bound} 
It therefore only remains to prove \( A\geq \Lambda \). 
Suppose by contradiction that $A < \Lambda$. 
We will show that this implies that there is a measure \( \mu \) 
for which 
\begin{equation}\label{contra}
\int \Phi d\mu <\Lambda
\end{equation}
which gives a contradiction. 
\subsubsection*{Construction of the measure \( \mu \)}
Notice first of all that, 
since \( n^{-1}A_{n}\to A \) for a.e. \( \omega \), we can
choose a set \( U \) of arbitrarily large measure on which this
convergence is uniform.  From this and the definition of \( A_{n} \), 
for every \( \epsilon > 0 \)
there exists \( N > 0 \) such that for all \( n\geq N \) there exists a measurable function \(
\omega \mapsto (x_n(\omega),
v_n(\omega)) \in SM  \) defined in \(U\), such that 
\begin{equation}\label{sub}
\frac{1}{n}\log
|D_{x_n(\omega)}\phi^{(n)}_{\omega} v_n(\omega)| = 
\frac{1}{n}A_{n}(\omega)  < A +\epsilon. 
\end{equation}
To see this, just consider the weakly measurable and closed valued set function 
\[ w\mapsto \{ (x,v)\in SM: \log |D_x\phi_w^{(n)}v|\, \mbox{is minimal}\} ,\]
 defined in \( U\), and choose any measurable selection  \( (x_n(\omega),
v_n(\omega))\) (for the existence of such a selection see, for instance, theorem 4.1 in \cite{W77}).
Then, for each \( \omega\in U \) and each \( n\geq 1 \) we define a 
probability measure 
\(
\sigma_{n}(\omega) = \delta_{x_{n}(\omega), v_{n}(\omega)}
\)
where \( \delta_{x,v} \) denotes the Dirac-delta measure at the point \( 
(x,v)\in SM \). 
We also let $G = \cup_{i=-\infty}^{\infty} 
\theta^i(U)$ (notice that ergodicity implies that $\mathbb{P}(G) = 1$)
and, for \( \omega\in G\setminus U \) define 
\(
\sigma_{n}(\omega) \equiv \delta_{x,v}
\)
for some arbitrary point \( (x,v)\in SM \) which can be chosen
independently of \( \omega \) or \( n \). 
We can now define 
probability measures 
$$
\mu_n(\omega) = \frac{1}{n}
\sum\limits_{i=0}^{n-1} ((\widehat{TF})^i \sigma_n)(\omega).
$$  
It is
easy to prove that the marginal of  $\mu_n$   on $ \Omega $
coincides with
  $ \mathbb{P}$ and it is  well known that $\{\mu_n \}
$ has a subsequence converging to an invariant measure $\mu \in
\mathcal M_{\mathbb{P}}(\widehat{TF}) $  (see Arnold
\cite{arnold98}). Without loss of generality, we suppose $
\lim\limits_{n \to \infty} \mu_n = \mu$. 

\subsubsection*{Contradiction}
It remains to prove \eqref{contra}, i.e. 
\(
\int_{\Omega \times SM} \Phi
d\mu < \Lambda,
\)
to get the desired contradiction. 
By the continuity of
\( \Phi \) we have 
\[ 
\int_{\Omega \times SM} \Phi
d\mu  = \lim_{n \to \infty} \int_{\Omega \times SM}
\Phi d\mu_n
\]
By the definition of \( \mu_{n} \) we have 
\[ 
\int_{\Omega \times SM}
\Phi d\mu_n = 
\frac{1}{n}\sum_{i=0}^{n-1}\int_{\Omega\times SM} \Phi \ d
((\widehat{TF})^{i}\sigma_{n})(\omega) d\mathbb P,
\]
and by the definition of \( \sigma_{n} \) and the fact that 
\( \mathbb P(G) =1 \), 
the right hand side above 
is equal to 
\[ 
\frac{1}{n}\int_{U}
\sum_{i=0}^{n-1} \Phi((\widehat{TF})^i(\omega,x_n(\omega),
v_n(\omega)))d\mathbb{P}+ \frac{1}{n}\int_{G\setminus U}
\sum_{i=0}^{n-1}
\Phi((\widehat{TF})^i(\omega,x, v))d\mathbb{P} 
\]
It is therefore sufficient to consider the limits of these two
integrals and show that their sum is strictly less than \( \Lambda \).
To bound the first part notice that 
\[ 
\sum_{i=0}^{n-1} \Phi((\widehat{TF})^i(\omega,x_n(\omega),
v_n(\omega))) = \log
|D_{x_n(\omega)}\phi^{(n)}_{\omega} v_n(\omega)| 
\]
and therefore, \eqref{sub} gives an upper bound of \( A+\varepsilon \).
For the second we have 
\begin{align*}
\int_{G\setminus U}
\sum_{i=0}^{n-1}
\Phi((\widehat{TF})^i(\omega,x, v))d\mathbb{P} 
& \leq \sum_{i=0}^{n-1}\int_{G\setminus U} 
|D\phi_{\theta^{i}(\omega)}| d\mathbb P 
\\ &= \sum_{i=0}^{n-1}\int_{\theta^{-i}(G\setminus U)}
|D\phi_{\omega}| d\mathbb P.
\end{align*}
Recall that \( |D\phi_{\omega}|=\max_{(x,v)}|D_{x}\phi_{\omega}(v)| \). 
 Since $\theta: \Omega \to \Omega $ is an invertible
 transformation preserving the ergodic measure $\mathbb{P}$, we have
 $\mathbb{P}(\theta^{-i}(G\setminus U)) = \mathbb{P}(G\setminus U)
 < \delta$. Thus, by the integrability condition on 
 \( |D\phi_{\omega}| \) we can choose \( U \) of sufficiently large
 measure so that 
$$\int_{\theta^{-i}(G\setminus U)}
|D\phi_{\omega}|  d\mathbb{P} <  \epsilon.$$  
Since \( \epsilon \) is arbitrary we get the desired contradiction. 
\end{proof}

\section{Tempered random variables}
Finally, it only remains to show that that ``constant'' \( C(\omega) \)
is a tempered random variable. To see this, notice first of all that
we can choose
\[ 
C(\omega) = \inf_{n \ge 1 } \left\{
e^{-\lambda n} \min_{(x,v) \in SM
 }|D_x\phi_{\omega}^{(n)}(v)| 
 \right\}
\]
Then we have
\begin{lemma}
    \( C(\omega) \) is a tempered random variable.
\end{lemma}
\begin{proof}
    We want to compare \( C(\theta
    \omega)\) to \( C(\omega) \). Let 
    \(
    D_{n}(\omega): = \min_{(x,v) \in SM
 }|D_x\phi_{\omega}^{(n)}(v)| .
    \)
  Then 
  \[ 
  \frac{C(\theta\omega)}{C(\omega)} = 
 \frac{\inf_{n \ge 1 } \left\{
 e^{-\lambda n} D_{n}(\theta\omega)
  \right\}}{\inf_{n \ge 1 } \left\{
 e^{-\lambda n} D_{n}(\omega)
  \right\}}
  = \frac{\min\{e^{-\lambda} D_{1}(\theta\omega), \inf_{n\geq
  2} \{e^{-\lambda n} D_{n}(\theta\omega)\}  \}}
  {\min\{e^{-\lambda} D_{1}(\omega), \inf_{n\geq
  2} \{e^{-\lambda n} D_{n}(\omega)\}  \}}.
  \]
 We consider two cases. Suppose first  that  
 \(
 C(\omega) = e^{-\lambda} D_{1}(\omega) \leq \inf_{n\geq
  2} \{e^{-\lambda n} D_{n}(\omega)\}
  \) 
  Then 
\[ 
\frac{C(\theta\omega)}{C(\omega)} \leq
\frac{D_{1}(\theta\omega)}{D_{1}(\omega)}.
\]
On the other hand, suppose that 
\(
C(\omega) = \inf_{n\geq
 2} \{e^{-\lambda n} D_{n}(\omega)\}\leq e^{-\lambda} D_{1}(\omega). 
\)
Then, keeping in mind that 
\( D_{n}(\omega) = D_{n-1}(\theta\omega)
D_{1}(\omega) \), we have, for any \( n\geq 2 \), 
\[ 
e^{-\lambda n}D_{n}(\omega) \geq e^{-\lambda
(n-1)}D_{n-1}(\theta\omega) e^{-\lambda} D_{1}(\omega) 
\geq C(\theta\omega) e^{-\lambda} D_{1}(\omega).
\]
Hence \( C(\omega) \geq  C(\theta\omega) e^{-\lambda} D_{1}(\omega)\) 
and so, 
combining the estimates in the two cases,  we have 
\[ 
\frac{C(\theta\omega)}{C(\omega)} \leq
\max\left\{\frac{D_{1}(\theta\omega)}{D_{1}(\omega)}, 
\frac{e^{\lambda}}{D_{1}(\omega)}
\right\}\leq \frac{\max\{ |D\phi_{\theta\omega}|, e^{\lambda} \}}
{D_{1}(\omega)}.
\]
Since 
\(
{1}/{D_{1}(\omega)} \leq |D\phi_{\omega}^{-1}|
\)
the integrability assumptions \eqref{intass} 
imply that
 \[
 \log^{+} \frac{C(\theta \omega )}{C(\omega)}
 \in L^1(\Omega,\mathcal{F}, \mathbb{P}).
 \]
 where $\log^+a = \max \{\log a, 0\}$.
The statement in the Lemma then follows by a
standard general result that any 
positive finite measurable function $g$ 
such that
$\log^+ \frac{g(\theta(\omega))}{g(\omega)}
\in L^1(\Omega,\mathcal{F}, \mathbb{P}) $
is tempered. For completeness we give a proof here. 
By the subadditive 
ergodic theorem the following
limit exists for a.e. $\omega $:
$$ \lim\limits_{n \to \infty} \frac1n \sum_{k =0 }^{n-1} \log
\frac{g \circ \theta^{k+1}}{g \circ \theta^k} = \lim\limits_{n \to
\infty} \frac1n \log \frac{g \circ \theta^n}{g} = 
\lim\limits_{n \to
\infty} \frac1n \log (g \circ \theta^n)=
h.
$$
The last equality follows from the fact that 
$\lim_{n \to \infty} n^{-1} \log g = 0 $. By the
definition of a tempered random variable, it is therefore 
sufficient to show
that \( h=0 \) for a.e. \( \omega \). 
For each fixed $\delta > 0 $ and using the invariance of the measure \( 
\mathbb P \) for \( \theta \), we have 
\begin{align*}
\lim\limits_{n \to \infty} \mathbb{P}(\{ \omega : \frac1n \log|g
\circ \theta^n(\omega)| \ge \delta \}) &=\lim\limits_{n \to
\infty} \mathbb{P}(\theta^{-n}g^{-1}(-e^{n\delta},e^{n\delta})^c)\\
&= \lim\limits_{n \to \infty} \mathbb{P}
(g^{-1}(-e^{n\delta},e^{n\delta})^c )= 0. 
\end{align*}
This means that the sequence of functions 
$\frac1n \log(g\circ \theta^n)$ converges to $0$
in measure and therefore some subsequence converges to 0 a.e. Since
we know from the above that the sequence actually converges a.e. this 
yields the result. 
\end{proof}

  \end{document}